\documentclass[a4paper,12pt]{amsart}
\usepackage{amssymb}
\usepackage{amscd}

\def\Bbb{\mathbb}
\def\Cal{\mathcal}

\def\endrk{\hbox{$|\!\!|\!\!|\!\!|\!\date!|\!\!|\!\!|$}}

\let\d=\delta
\let\e=\varepsilon

\let\i=\iota

\def\M{{\Cal{M}}}

\def\LOT{{\rm LOT}}

\newcommand{\IT}[1]{{\rm(}{\it{\!#1}}{\rm)}}

\newcommand{\newc}{\newcommand}

\let\ccdot\cdot
\def\cdot{\hbox to 2.5pt{\hss$\ccdot$\hss}}

\newcommand{\om}{\omega}
\renewcommand{\phi}{\varphi}

\newcommand{\si}{\sigma}

\newcommand{\Om}{\Omega}

\newc{\aI}{\mbox{\boldmath{$ I$}}}
\newc{\aR}{\mbox{\boldmath{$ R$}}}
\newc{\aDeR}{\mbox{\boldmath{$ U$}}_B{}^P{}_C{}^Q}
\newc{\al}{\mbox{\boldmath$ \Delta$}}             
\newc{\nda}{\mbox{\boldmath$ \nabla$}}
\newc{\ad}{\mbox{\boldmath$ d$}}
\newc{\da}{\mbox{\boldmath$ \delta$}}
\newc{\aK}{\mbox{\boldmath{$ K$}}}
\newc{\aL}{\mbox{\boldmath{$ L$}}}

\newtheorem{theorem}{Theorem}[section]
\newtheorem{lemma}[theorem]{Lemma}
\newtheorem{proposition}[theorem]{Proposition}
\newtheorem{corollary}[theorem]{Corollary}

\newcommand{\cB}{{\Cal B}}

\newcommand{\ce}{{\Cal E}}

\newcommand{\cq}{{\Cal Q}}

\newcommand{\nd}{\nabla}

\newcommand{\Rho}{{\mbox{\sf P}}}
\newcommand{\Up}{\Upsilon}

\newcommand{\Aut}{\operatorname{Aut}}
\newcommand{\End}{\operatorname{End}}

\newcommand{\Ric}{\operatorname{Ric}}

\renewcommand{\S}{{\Bbb S}}

\newcommand{\R}{\mathbb{R}}

\newcommand{\nds}{D}

\newcommand{\Mtw}{{M^{\nds}}}

\newcommand{\dtw}{{d^{\nds}}}
\newcommand{\twd}{{\d^{\nds}}}

\newcommand{\dtwpm}{{d^{\nds}_\pm}}

\newcommand{\twdpm}{{\d^{\nds}_\pm}}

\newcommand{\mstar}{\mbox{\large $ \star$}}

\newcommand{\bD}{\Bbb D}

\let\i=\iota

\newcommand{\nn}[1]{(\ref{#1})}

\newcommand{\D}{\mbox{\boldmath{$ D$}}}

\newcommand{\bT}{\Bbb T}

\newcommand{\h}{\mbox{\boldmath{$ h$}}}
\newcommand{\bg}{\mbox{\boldmath{$ g$}}}

\def\cD{{\Cal D}}

\newcommand{\V}{{\mbox{\sf P}}}                   
\newcommand{\J}{{\mbox{\sf J}}}

\newc{\strutdd}{\rule{0mm}{5mm}}

\usepackage{ifthen}

\newc{\tensor}[1]{#1}

\newc{\Mvariable}[1]{\mbox{#1}}

\newc{\down}[1]{{}_{
\ifthenelse{\equal{#1}{;}}{|}{#1}}}

\newc{\up}[1]{{}^{#1}}
\newc{\C}{C}

\newc{\JulyStrut}{\rule{0mm}{6mm}}
\newc{\midtenPan}{\mbox{\sf S}}
\newc{\midten}{\mbox{\sf T}}
\newc{\midtenEi}{\mbox{\sf U}}
\newc{\ATen}{\mbox{\sf E}}
\newc{\BTen}{\mbox{\sf F}}
\newc{\CTen}{\mbox{\sf G}}

\def\sideremark#1{\ifvmode\leavevmode\fi\vadjust{\vbox to0pt{\vss
 \hbox to 0pt{\hskip\hsize\hskip1em
 \vbox{\hsize3cm\tiny\raggedright\pretolerance10000
 \noindent #1\hfill}\hss}\vbox to8pt{\vfil}\vss}}}%

\input epsf
\input xy
\xyoption{all}

\newtheorem{definition}[theorem]{Definition}

\newtheorem{remark}[theorem]{Remark}

\newcommand{\bd}{\begin{definition}}
\newcommand{\ed}{\end{definition}}

\newcommand{\br}{\begin{remark}}
\newcommand{\er}{\end{remark}}

\newcommand{\bt}{\begin{tabular}}
\newcommand{\et}{\end{tabular}}

\def\cdot{\hbox to 2.5pt{\hss$\ccdot$\hss}}

\usepackage{amsmath,amsfonts,amssymb}

\begin{document}
\renewcommand{\today}{}
\title{Yang-Mills detour complexes and conformal geometry}
\author{A.\ Rod Gover, Petr Somberg and Vladim\'{\i}r Sou\v cek}

\address{ARG: Department of Mathematics\\
  The University of Auckland\\
  Private Bag 92019\\
  Auckland 1\\
  New Zealand} \email{gover@math.auckland.ac.nz}
\address{PS and VS: Mathematical Institute\\
Faculty of Mathematics and Physics\\
Charles University,
Sokolovsk\'a 83, 186 75 Praha\\
Czech Republic}\email{somberg@karlin.mff.cuni.cz,soucek@karlin.mff.cuni.cz}

\vspace{10pt}

\thanks{ARG would like to thank
the Royal Society of
New Zealand for support via Marsden Grants no.\ 02-UOA-108 and 06-UOA-029. 
Part of the work was prepared during a visit of ARG
supported by the E.\ \v Cech Center. PS and VS acknowledge the support of
the grant GA CR 201/05/2117 and the grant MSM 021620839}

\renewcommand{\arraystretch}{1}
\maketitle
\renewcommand{\arraystretch}{1.5}

\pagestyle{myheadings}
\markboth{Gover, Somberg, Sou\v cek}{Yang-Mills detour complexes}

\begin{abstract} 
Working over a pseudo-Riemannian manifold, for each vector bundle with
connection we construct a sequence of three differential operators
which is a complex (termed a Yang-Mills detour complex) if and only if
the connection satisfies the full Yang-Mills equations. A special case
is a complex controlling the deformation theory of Yang-Mills
connections.  In the case of Riemannian signature the complex is
elliptic. If the connection respects a metric on the 
bundle then the complex is formally self-adjoint. In dimension 4 the
complex is conformally invariant and generalises, to the full
Yang-Mills setting, the composition of (two operator) Yang-Mills
complexes for (anti-)self-dual Yang-Mills connections. Via a prolonged
system and tractor connection a diagram of differential operators is
constructed which, when commutative, generates differential complexes
of natural operators from the Yang-Mills detour complex. In dimension
4 this construction is conformally invariant and is used to yield two
new sequences of conformal operators which are complexes if and only
if the Bach tensor vanishes everywhere. In Riemannian signature these
complexes are elliptic.  In one case the first operator is the twistor
operator and in the other sequence it is the operator for Einstein
scales. The sequences are detour sequences associated to certain
Bernstein-Gelfand-Gelfand sequences.
\end{abstract}

\section{Introduction}

In the study of Riemannian and pseudo-Riemannian geometry it is often
valuable to use differential operators with good conformal behaviour.
In the Riemannian setting, elliptic differential operators are
particularly important. For example the conformal Laplacian controls
the conformal variation of the scalar curvature. This was exploited
heavily in the solution by Schoen, Aubin, Trudinger, and Yamabe (see
\cite{L&P}) of the ``Yamabe Problem'' of finding, via conformal
rescaling, constant scalar curvature metrics on compact manifolds.
Related curvature prescription problems and techniques have exploited
the higher order conformal Laplacians of Paneitz, Graham et al.\
\cite{GJMS,Brsharp,aliceC-EMS}. These operators on functions (or
really densities) also find a natural place in the recent developments
\cite{FGrnew,GrZ} concerning the asymptotics and scattering theory of
the conformally compact Poincar\'e-Einstein metric of Fefferman-Graham
\cite{FGrast}.

On many tensor and spinor fields there is no conformally invariant
elliptic operator (taking values in an irreducible bundle); this
follows from the classification of conformally invariant differential
operators on the sphere \cite{BoeColl,EastSlo}. This classification is
based on the structure of generalised Verma modules and from this it
follows that often the analogue, or replacement, for a conformal
elliptic operator on the sphere is an elliptic complex of conformally
invariant differential operators.  However the situation is
complicated for conformally curved structures. The requirement that a
sequence of differential operators be both conformally invariant and
form a complex is severe. On the other hand when such complexes exist
they can be expected to play a serious role in treating the underlying
structure.  This idea is already well-established in the setting of
self-dual 4-manifolds \cite{AHS,Donaldson}. On fully conformally
curved $n$-manifolds, with $n$ even, there is a class of elliptic
conformal complexes on differential forms \cite{BrGodeRham}. Each of
these is different to the de Rham complex, and these complexes
generalise the conformally invariant operator of \cite{GJMS}, with
leading term $\Delta^{n/2}$.  Another class of complexes is
based around the (Fefferman-Graham) obstruction tensor \cite{FGrast}.
This is a natural conformal 2-tensor that generalises, to higher even
dimensions, the Bach tensor in dimension 4.  It turns out that the
formal deformations of obstruction-flat manifolds are controlled by a
sequence of conformal operators, which form an elliptic complex if and
only if the structure is obstruction-flat \cite{defdet}.
Unfortunately there is no obvious way to generalise either the
construction in \cite{BrGodeRham}, or that in \cite{defdet}.

For 4-manifolds we construct
here two conformal differential sequences 
which are (formally
self-adjoint) complexes if and only if the (conformally invariant) Bach-tensor \cite{Bach} vanishes
everywhere. This condition is weaker
than self-duality. In fact conformally Einstein manifolds are also
Bach-flat and there are structures which are Bach-flat and neither
conformally-Einstein nor half-flat
\cite{GoLeit}. Writing ${\bf T}: \S \to {\rm Tw} $ for the usual twistor
operator on Dirac spinors (as in e.g. \cite{BFGK}), in Theorem
\ref{twcase} we obtain a differential complex
$$
 \S  \stackrel{\bf T}{\to} {\rm Tw} 
\stackrel{M^\Sigma}{\longrightarrow} {\rm Tw}
\stackrel{{\bf T}^*}{\to} {\S} ,
$$ where $M^\Sigma$ is a third order Rarita-Schwinger type 
operator.  
On the other hand in Theorem \ref{paracase} we construct
$$
\ce^0\stackrel{P}{\to} \ce^{1,1} 
\stackrel{M^{\bT}}{\longrightarrow} \ce^{1,1}
\stackrel{P^*}{\to} \ce^0
$$ where $M^{\bT}$ is a second order conformal operator, similar in
form to the operator which controls deformations of Einstein
structures (see \cite{Besse} and references therein), while $P$ is a
curvature modification of the trace-free covariant Hessian.
Non-vanishing solutions of $P$ give conformal factors $\si$ so that
$\si^{-2}g$ is Einstein (see \cite{BEGo}); we show via the second
sequence that the Bach tensor obstructs solutions.  If the manifold is
Riemannian then both of the complexes are elliptic.  We have been
intentionally explicit in treating these constructions, as it seems
these complexes should play a fundamental role in conformal and
Riemannian geometry.  In the compact and Riemannian-signature setting
the ellipticity implies that the complexes have finite dimensional
cohomology spaces.  In both cases the interpretation of the $0^{\rm
  th}$-cohomology is well-known but as far as we know the first
cohomology is a new global conformal invariant of Bach-flat
structures.  Such conformal elliptic complexes have the scope to yield
further geometric information through their detour torsion invariants
\cite{Tsrni}.

In fact the Theorems \ref{paracase} and \ref{twcase} construct short
detour complexes in all dimensions $n\geq 3$ and $n\geq 4$
respectively. These complexes are conformally invariant only in
dimension 4, but by construction have a simple conformal behaviour and
may well be of interest for physics in the Lorentzian setting.

The route to the constructions and results mentioned above is really
one of the main points of the article. We believe that it lays
foundations for an eventual general treatment of a large class of related 
complexes,
and also many of the results should be of independent
interest.  The simplest example of a detour complex is the Maxwell
detour complex
\begin{equation}\label{Maxdet}
\ce^0 \stackrel{d}{\to} \ce^1 \stackrel{\d d}{\to} 
\ce^1\stackrel{\d}{\to}\ce^0 .
\end{equation} 
For each vector bundle $V$ and connection $D$ we construct, in section
\ref{gencon}, a curvature adjusted twisting of this complex with the
property that it is again a complex if and only if the connection $D$
is a (pure) Yang-Mills connection, see Theorem \ref{twistthm}.  In
dimension 4 the resulting complexes are conformal.  In sections
\ref{var} (in particular Theorem \ref{def}) we recover a class of
these complexes by considering deformations of Yang-Mills connections.
We show in Proposition \ref{agree} of Section \ref{hflat} that, in
dimension 4, the Yang-Mills detour complex generalises the composition
of subcomplexes of Yang-Mills complexes arising from (anti-)self-dual
connections.  The next main item is a rather general construction, see
diagram [D] in Section \ref{detdiag}.  This enables the Yang-Mills
detour complex to be ``translated'' to yield new complexes. Broadly
the motivational idea is this. If one has an overdetermined
differential operator (of finite type) $\cB^0\to \cB^1$ then one may
sometimes obtain a corresponding invariant connection on a prolonged
system \cite{GoSil}. If the latter satisfies the Yang-Mills equations
and, say, preserves a metric on the prolonged system then the
Yang-Mills detour complex on the prolonged system descends and extends
$\cB^0\to \cB^1$ to a complex.  In reality this is an
over-simplification, but it contains the germ of the main idea.

The first author would like to thank Helga Baum, Kengo Hirachi,
Paul-Andi Nagy and Andrew Waldron for illuminating
discussions.

\section{Background: conformal geometry}\label{back}

 Recall that a {\em conformal structure\/} of signature $(p,q)$ on $M$
is a smooth ray subbundle $\cq\subset S^2T^*M$ whose fibre over $x$
consists of conformally related signature-$(p,q)$ metrics at the point
$x$ (and $S^2 T^*M$ is the symmetric part of $\otimes^2 T^*M$). Sections of $\cq$ are metrics $g$ on $M$. So we may equivalently
view the conformal structure as the equivalence class $[g]$ of these
conformally related metrics.  The principal bundle $\pi:\cq\to M$ has
structure group $\Bbb R_+$, and so each representation ${\Bbb R}_+ \ni
x\mapsto x^{-w/2}\in {\rm End}(\Bbb R)$ induces a natural line bundle
on $ (M,[g])$ that we term the conformal density bundle $E[w]$. We
shall write $ \ce[w]$ for the space of sections of this bundle and
$\bg$ denotes the {\em conformal metric}, that is the tautological
section of $S^2T^*M\otimes E[2]$ determined by the conformal
structure. On conformal manifolds this will be used to identify $TM$
with $T^*M[2]$. Note $E[w]$ is trivialised by a choice of metric $g$
from the conformal class, and we write $\nd$ for the connection
corresponding to this trivialisation (and term this the Levi-Civita
connection on $E[w]$).  It follows that (the coupled) $ \nd_a$
preserves the conformal metric.

In dimensions $n\geq 3$ the Riemannian
curvature can be decomposed into the totally trace-free Weyl curvature
$C_{abcd}$ and a remaining part described by the symmetric {\em
Schouten tensor} $\Rho_{ab}$, according to $
R_{abcd}=C_{abcd}+2\bg_{c[a}\Rho_{b]d}+2\bg_{d[b}\Rho_{a]c}, $ where
$[\cdots]$ indicates antisymmetrisation over the enclosed indices.
The Schouten tensor is a trace modification of the Ricci tensor
$\Ric_{ab}$ and vice versa: $\Ric_{ab}=(n-2)\Rho_{ab}+\J\bg_{ab}$,
where we write $ \J$ for the trace $ \V_a{}^{a}$ of $ \V$.  The {\em
Cotton tensor} and {\em Bach tensor} are defined by, respectively, 
\begin{equation} 
A_{abc}:=2\nabla_{[b}\Rho_{c]a} \quad \mbox{and} \quad B_{ab}:=\nabla^c
A_{acb}+\Rho^{dc}C_{dacb}.\label{bach} 
\end{equation}

Under a {\em conformal transformation} we replace a choice of metric $
g$ by the metric $ \hat{g}=e^{2\om} g$, where $\omega$ is a smooth
function. Explicit formulae for the corresponding transformation of
the Levi-Civita connection and its curvatures are given in e.g.\
\cite{BEGo,GoPetLap}. We recall that, in particular, the Weyl
curvature is conformally invariant $\widehat{C}_{abcd}=C_{abcd}$.  In
dimension 4 $B_{ab}$ is conformally invariant.

We will write $\ce^k[w]$ for the sections of the tensor product
$E^k[w]:=\wedge^kT^*M \otimes E[w]$.  
On conformal manifolds we use the notation $\ce_k$ to mean
the space of sections of $E_k:=\wedge^kT^*M \otimes E[2k-n]$.  This
notation (following \cite{BrGodeRham}) is suggested by the duality
between the section spaces $\ce^k$ and $\ce_k $; compactly supported
sections pair globally by contraction and integration.  For any vector
bundle $V$, $\ce^k(V)$ is the space of smooth sections of
$E^k(V):=\wedge^kT^*M \otimes V$, while $\ce_k(V)$ means the space of sections
of $E_k(V):=\wedge^kT^*M \otimes E[2k-n] \otimes V$. When a metric from the
conformal class is fixed, these spaces will be identified.

In conformal geometry the de Rham complex is a prototype for a
class of sequences of bundles and conformally invariant differential
operators, each of the form 
$$
\cB^0\to \cB^1\to \cdots \to \cB^n~, 
$$ where the vector bundles $ \cB^i$ are irreducible tensor-spinor
bundles. On the $n$-sphere there is one such complex for each irreducible
module $\mathbb{V}$ for the group $G=SO(n+1,1)$ of conformal
motions, the space of solutions of the first (overdetermined)
conformal operator $\cB^0 \to \cB^1$ is isomorphic to $ \mathbb{V}$, and the 
sequence
gives a resolution of this space viewed as a sheaf.  These are the
conformal cases of the (generalised) Bernstein-Gelfand-Gelfand (BGG)
sequences, a class of sequences of differential operators that exist
on any parabolic geometry \cite{BoeColl,CSSannals}.  As well as
the operators $D_i:\cB^i\to \cB^{i+1}$ of the BGG sequence, in even
dimensions there are conformally invariant ``long operators'' $
L_k:\cB^k\to \cB^{n-k}$ for $ k=1,\cdots ,\small{n/2-1}$ \cite{BoeColl}.
Thus there are sequences of the form
$$
\cB^0\stackrel{D_0}{\to}\cB^1 \stackrel{D_1}{\to}\cdots 
\stackrel{D_{k-1}}{\to} \cB^k \stackrel{L_k}{\to} \cB^{n-k} 
\stackrel{D_{n-k}}{\to}\cdots \stackrel{D_{n-1}}{\to} \cB^n~.
$$ and, following \cite{BrGodeRham,Tsrni}, we term these detour sequences since, in
comparison to the BGG sequence, the long operator here bypasses the middle of the BGG sequence. Once again from the classification it follows
that these detour sequences are in fact complexes in the
case that the structure is conformally flat. 
  The dimension 4 conformal complexes, constructed in 
Theorems \ref{paracase} and \ref{twcase} below,
 are detour sequences of this form with $k=1$.

\section{Yang-Mills detour complexes} \label{YMdets}

\subsection{The general construction} \label{gencon}

We work over a pseudo-Riemannian $n$-manifold $(M,g)$ of signature
$(p,q)$ ($n\geq 2$).  Let $V$ denote a vector bundle with a connection
$\nds$. We denote by $F$ the curvature of $\nds$.  We also write
$\nds$ for the induced connection on the dual bundle $V^*$.  We write
$d^\nds$ for the connection-coupled exterior derivative operator
$\dtw:\ce^k(V)\to \ce^{k+1}(V)$. Of course we could equally consider
$\dtw:\ce^k(V^*)\to \ce^{k+1}(V^*)$, and for the formal adjoint of
this we write $\d^\nds:\ce_{k+1}(V)\to\ce_k(V)$.

Let us write $F\cdot$ for the action of
the curvature on the twisted 1-forms, $F\cdot :\ce^1(V)\to \ce_1(V)$
given by 
$ 
(F\cdot \phi)_a:= F_a{}^b\phi_b
$ where we have indicated the abstract form indices explicitly, whereas
the standard $\End(V)$ action of the curvature on the $V$-valued
1-form is implicit. Using this we construct a
differential operator
$$
M^{\nds}:\ce^1(V)\to \ce_1(V)
$$
by 
$
M^{\nds}\phi=\d^\nds d^\nds\phi-F\cdot \phi .
$

The operator $M^\nds$ has the property that its composition with
$d^\nds$ is given simply by an algebraic action of the ``Yang-Mills
current'' $\delta^\nds F$ on the bundle $V$, as follows.
\begin{lemma}\label{algact}
The composition of $\Mtw: \ce^1(V)\to \ce_1(V)$ with $\dtw:
\ce^{0}(V)\to \ce^1(V)$ is given by the exterior action of
$\d^{\nds}F$, as an $\End(V)$-valued 1-form:
$$
\Mtw \dtw = \e(\twd F).
$$

The composition of $\twd: \ce_1(V)\to \ce_0(V)$ with 
$\Mtw: \ce^1(V)\to \ce_1(V)$ is given by the interior action of
$-\d^{\nds}F$,  
 as an $E^1$-valued endomorphism of $E^1(V)$:
$$
\twd \Mtw  = - \i(\twd F).
$$

In these expressions the interior multiplication (indicated by
$\i(\cdot)$) and the exterior multiplication (indicated $\e(\cdot)$)
refers to the form index of $\twd F$.
\end{lemma}
\noindent{\bf Proof:} For the connection $\nds$ coupled with the
Levi-Civita connection $\nd$, let us also write $\nds$. Then, again using the 
notation where we exhibit abstract tensor indices but
suppress indices for the bundle $V$,  a formula for $\Mtw$ on a twisted 1-form
$\Psi_a$ is 
$$
(\Mtw \Psi)_b =-\nds^a\nds_{a}\Psi_{b}+\nds^a\nds_{b}\Psi_{a}-F_{b}{}^a\Psi_a,
$$ since the Levi-Civita connection is torsion-free. On the other hand
for $\Phi\in\ce^0(V)$, $(\dtw \Phi)_a=\nds_a \Phi$. Thus
$$
\begin{aligned}
(\Mtw \dtw \Phi)_b 
&=\nds^a(\nds_{b}\nds_{a}\Phi- \nds_{a}\nds_{b}\Phi)-F_{b}{}^a\nds_a\Phi\\
&= \nds^a F_{ba}\Phi - F_{ba}\nds^a\Phi\\
&=\big(\e(\twd F)\Phi\big)_b ~.
\end{aligned}
$$

\vspace{2mm}

\noindent By a similar calculation (or using the above on $(V^*,D)$
and taking formal adjoints) we obtain,
$$
\twd \Mtw \Psi  =- \i(\twd F)\Psi,
$$ for $\Psi\in \ce_1(V)$. (Note that $ \i(\twd F)\Psi=-(\nds^b
F_{b}{}^a) \Psi_a $.)  \quad $\blacksquare$\\ 
\noindent{\bf Remark:} 
 Note that to simplify the punctuation in calculations, we often view sections 
of vector bundles as order 0 operators.  Thus for example 
$\nds^a F_{ba}\Phi$ has the same meaning as $\nds^a (F_{ba}\Phi) $.
\quad \endrk 

\vspace{2mm}

If the connection $\nds$ is orthogonal or unitary for some inner
product or Hermitian form on $V$ (then $V$ may be identified with
$V^*$ and) the algebraic action $F\cdot :\ce^1(V)\to \ce_1(V)$ is
easily verified to be formally self-adjoint and so, in this case,
$\Mtw$ is formally self-adjoint.  
From these observations, and the Lemma \ref{algact}, we have
the following.
\begin{theorem}\label{twistthm}
The sequence of operators,
\begin{equation}\label{detseq}
\ce^0(V)\stackrel{d^\nds}{\to} \ce^1(V)\stackrel{M^\nds}{\to} 
\ce_1(V)\stackrel{\d^\nds}{\to}\ce_0(V)
\end{equation}
is a complex if and only if 
the curvature $F$ of the connection $\nds$ satisfies the (pure) Yang-Mills
equation
$$
\d^\nds F=0.
$$ In addition:\\ \IT{i} If $\nds$ is an orthogonal or unitary
connection then the sequence is formally self-adjoint.\\ \IT{ii} In
Riemannian signature the sequence is elliptic.\\
\IT{iii} In dimension 4 the sequence \nn{detseq} is conformally
invariant.
\end{theorem}
\noindent{\bf Proof:} It remains to show \IT{ii} and \IT{iii}. For
\IT{ii} we need that the symbol sequence is exact. This sequence is
simply a tensor product twisting by $V$ of the symbol sequence of the
Maxwell detour complex \nn{Maxdet} and so it is sufficient to check
that case. But that case is an easy consequence of the algebraic Hodge
decomposition on an inner product space.

The conformally well-defined formal
adjoint of the exterior derivative $d:\ce^k\to \ce^{k+1}$ acts
$$
\d:\ce_{k+1}\to \ce_k .
$$ 
(cf.\ e.g.\ \cite{BrGodeRham}.) Note that in even dimensions on middle order
forms we have $\ce^{n/2}= \ce_{n/2}$ and so $\d:\ce^{n/2}\to
\ce_{n/2-1}$ is conformally invariant. The invariance persists if we twist
by a connection $\nds$,  and so from 
the definition of $M^D$ we have the result.  \quad $\blacksquare$

\vspace{1mm}

For a given connection $\nds$ on a vector bundle $V$, such that $\twd
F=0$, we will term the complex \nn{detseq} of Theorem \ref{twistthm} the
(corresponding) {\em Yang-Mills detour complex}.

If $\nds$ is a Yang-Mills connection on a vector bundle $V$, then the
dual connection on $V^*$ and the tensor product connection on any
tensor power of these are also Yang-Mills. One might alternatively
work with principal connections. If $\omega$ is a Yang-Mills
connection on a principal bundle ${\Cal P}$ with structure group $G$,
then we obtain a complex \nn{detseq} for every finite dimensional
representation of $G$.

\subsection{A variational construction of the deformation detour}
\label{var}
\newcommand{\A}{\dot{A}}

Returning to the general situation that began section \ref{gencon},
let $V$ denote a vector bundle with a connection $\nds$ and denote by
$F$ the curvature of $\nds$. Consider now a smoothly parametrised family of
connections $\nds^t$ (on $V$) given, on a section $v\in \ce^0(V)$, by 
\begin{equation}\label{vari}
\nds^t_a v = \nds_a v + A^t_a v
\end{equation}
where for each $t\in {\Bbb R}$, $A^t\in \ce^1({\End} V)$ and $A^0=0$. With 
$F^t$ denoting the curvature of $\nds^t$, we have 
$$
F^t_{ab}= F_{ab}+ \nds_a A^t_b-\nds_b A^t_a + [A^t_a , A^t_b] ,
$$ where,  once again, we write $\nds$ also to mean the connection on
$V$ coupled with the Levi-Civita connection.  It follows that 
the derivative of
$F^t$ at $\nds=\nds^0$ is
$$ \dot{F}_{ab}=\nds_a \A_b -\nds_b \A_a \quad \mbox{ where } \quad \A_a:=
\frac{d}{dt} A^t_a|_{t=0}~,
$$
that is $\dot{F}=d^\nds \A $. Now we calculate the derivative, at $\nds$, of 
$\delta^{\nds^t} F^t$. We have 
$$
\begin{aligned}
\frac{d}{dt} g^{ab}\nds^t_a F^t_{bc}|_{t=0} &=\nds^b \dot{F}_{bc}+[\A^b,F_{bc}] \\
&= \nds^b (\nds_b \A_c -\nds_c \A_b ) + [F_c{}^b,\A_b]~,
\end{aligned}
$$
where $\A$ acts on $F$ and vice versa by the obvious composition of bundle 
endomorphisms. Note that, since the 1-form $\A$ has values in ${\End} V$, 
the last term here is $F\cdot \A$.
Multiplying the display by $-1$ gives 
\begin{equation}\label{currentder}
\frac{d}{dt} \delta^{\nds^t} F^t |_{t=0}= M^{\nds} \A . 
\end{equation}
So we have, in particular, the following outcome.
\begin{lemma} 
If $\nds$ is a Yang-Mills connection then the infinitesimal
deformation $\A$ of $\nds$ is through Yang-Mills connections if and
only if $M^{\nds} \A =0$.
\end{lemma}

In the vector bundle picture, a so-called gauge transformation arises
locally by acting on $V$ by a section $u$ of the fibre bundle
$\Aut (V)$ of invertible elements in $\End (V)$.  
From the Leibniz rule
for $\nds$ (viewed as a connection on the tensor powers of  $V$ and $V^*$)
it follows immediately that
this pulls back to a transformation 
$$
\nds_a  \mapsto \nds_a +u^{-1} \nds_a u ,
$$ 
of the  connection, and 
whence 
\begin{equation}\label{curvt}
F_{ab} \mapsto u^{-1}F_{ab} u, \quad \mbox{ and } \quad \nds^a F_{ab}\mapsto u^{-1}(\nds^a F_{ab}) u~. 
\end{equation}
Thus if $u_s$ is a smoothly parametrised family of such
transformations with $u_0=id_V$ and derivative 
$$ 
\frac{d}{ds} u_s|_{s=0}=\dot{u}\in \ce^0(\End (V))
$$
then we obtain that the infinitesimal variation of $\nds^s$ is exactly $d^\nds \dot{u}$:
\begin{equation}\label{gaugeder}
\dot{\nds}_a=\nds_a \dot{u} ~.
\end{equation}
So from this and \nn{currentder} we have
$$
\frac{d}{ds} \delta^{\nds^s} F^s |_{s=0}= M^{\nds} d^\nds \dot{u} .
$$
On the other hand from \nn{curvt} and \nn{gaugeder} we get 
$$
\frac{d}{ds} \delta^{\nds^s} F^s |_{s=0}= (\delta^{\nds} F)\dot{u} -\dot{u} \delta^{\nds} F .
$$
Putting the last two results together brings us to 
$$
M^{\nds} d^\nds \dot{u} = \e( \delta^{\nds} F^{\End (V)}) \dot{u} 
$$ where $F^{\End (V)}$ is the curvature of $\nds$ viewed as a connection on
$\End (V)$ (so e.g.\ $F^{\End (V)} \dot{u}=[F,\dot{u}]$). This agrees precisely
with the specialisation of Lemma \ref{algact} to $\End (V)$ equipped
with the connection induced from $\nds$ on $V$.
In particular if $\nds $ is a Yang-Mills connection then so is the connection 
on $\End (V)$.
Since $\End (V)$ carries the non-degenerate symmetric pairing $(U,W)= Tr(UW)$
and this is preserved by
$\nds$, then it follows from \nn{gaugeder} that $M^\nds$ is formally
self-adjoint with respect to the global pairing obtained by integrating $(~,~)$. (The point is that the Yang-Mills equations are the Euler-Lagrange equations for the Lagrangian density $Tr(F^{ab}F_{ab})$. So  
 \nn{gaugeder} shows that $M^\nds$ is the second
variation of an action.  By interchanging orders of variation
one obtains the symmetry.) Thus from $M^\nds d^\nds$ 
we also have $\d^\nds M^\nds=0$ and the following result. 
\begin{theorem}\label{def}  For a vector bundle $V$, with Yang-Mills 
connection $\nds$, the (formal)
deformation detour complex
\begin{equation}\label{firstcx}
\ce^0(\End (V))\stackrel{d^\nds}{\to} \ce^1(\End (V))\stackrel{M^\nds}{\to} 
\ce_1(\End (V))\stackrel{\d^\nds}{\to}\ce_0(\End (V))
\end{equation}
is formally self-adjoint. Its first cohomology $H^1( \End (V),D) $ is the formal
tangent space at $\nds$ to the moduli space of Yang-Mills connections on $V$. 
 \end{theorem}
 \noindent It follows from a general deformation theory that 
 the complex \nn{firstcx}
 controls the full formal deformation theory of the Yang-Mills equations.

\subsection{Examples: (pseudo-)Riemannian manifolds with harmonic curvature}
On a pseudo-Riemannian (spin) manifold we write $\nd$ for the
Levi-Civita connection and $R$ for its curvature, the Riemannian
curvature tensor.  Riemannian structures satisfying $\d^{\nd}R=0$ are
said to have harmonic curvature. Einstein manifolds, for example, are
harmonic in this sense. There is a rich theory of harmonic manifolds,
see \cite{Besse} and references therein.

If $\d^{\nd}R=0$ then, from Theorem \ref{twistthm}, we get a detour
complex \nn{detseq} for $V$ any tensor (spin) bundle. For example if 
$T M$ is the tangent bundle then we have 
$$
M^\nd:\ce^1(TM)\to\ce_1(TM) \quad \mbox{by} \quad S_b{}^c\mapsto 
-2\nd^a\nd_{[a } S_{b]}{}^c- R_{ba}{}^c{}_dS^{ad} .
$$ This annihilates the covariant derivative of any tangent vector
field.

\subsection{Half-flat connections}\label{hflat}
In the setting of conformal (or pseudo-Riemannian)
4-manifolds, we observe here that
when a vector bundle connection $\nds$ is half-flat then there is very
simple interpretation of the Yang-Mills detour complex. First
we review, in our current notation, some relevant (well-known) background. 

Recall that on a conformal 4-manifold $\M$ of signature $(p,q)$ we
have $\mstar\mstar=(-1)^{k(4-k)+q}$ on $k$-forms. In the case of
Minkowskian signature let us write $E^2_{\pm}$ for the $\pm
i$-eigenspaces of $\mstar$. In the other signatures $E^2_{\pm}$ means
the $\pm 1$ eigenspaces of $\mstar$. In any case, since $\mstar$ is a
symmetric endomorphism of $E^2$, the decomposition of $E^2$ into
$E^2_+\oplus E^2_-$ is orthogonal.
Viewing the curvature $F$ (of $\nds$ on $V$) as a twisted 2-form,
recall that the curvature, or the connection, is said to be {\em self-dual}
(respectively {\em anti-self-dual}) if the component of $F$ in
$\ce^2_-(\End (V))$ (respectively in $\ce^2_+(\End (V))$) is zero, $F_-=0$
(respectively $F_+=0$).  So if a connection $\nds$ is half-flat, in
this sense, then $\twd F$ is a multiple of $\star d^\nds F$. But this
vanishes by the differential Bianchi identity for $F$. So $\twd F=0$
for connections which are either self-dual or anti-self-dual and each
case gives a special setting where the sequence \nn{detseq} is a
complex.

Let us write $\dtwpm$ for the compositions given by
$\dtw:\ce^1(U)\to \ce^2(U)$ followed by the projections $\ce^2(U)\to
\ce^2_\pm(U)$, where $U$ means either the bundle $V$ or its dual $V^*$. 
Thus by construction the operators $\dtwpm:\ce^1(U)\to
\ce^2_\pm(U)$ are conformally invariant. We write
$\twdpm:\ce^2_\pm(V)\to\ce_1(V) $  for the operators formally adjoint to 
$\dtwpm:\ce^1(V^*)\to
\ce^2_\pm(V^*)$. By construction these also 
are conformally invariant.
Now on $\Phi\in \ce^0(V)$ we have $\dtw\dtw \Phi=F\Phi$. The
projection of this into $\ce^{2}_\pm$ vanishes for all $\Phi$ if and
only if $F_{\pm}=0$.  By a similar observation for the composition
$\dtw \dtw$ on $\ce^0(V^*)$, and then taking formal adjoints, we see
that we have the situation in the following proposition. These results
are well-known. 
\begin{proposition}\label{subcxs} 
The sequences 
$$ \ce^0(V)\stackrel{\dtw}{\longrightarrow}
\ce^1(V)\stackrel{d^{\nds}_+}{\longrightarrow}\ce^2_{+}(V)
\quad {\rm and} \quad 
\ce^2_{+}(V) \stackrel{\d^\nds_+}{\longrightarrow} 
\ce_1(V) \stackrel{\d^\nds}{\longrightarrow} \ce_0(V)
$$ 
are complexes if and only if $F_{+}=0$. Similarly 
the sequences 
$$ \ce^0(V)\stackrel{\dtw}{\longrightarrow}
\ce^1(V)\stackrel{d^{\nds}_-}{\longrightarrow}\ce^2_{-}(V)
\quad {\rm and} \quad 
\ce^2_{-}(V) \stackrel{\d^\nds_-}{\longrightarrow} 
\ce_1(V) \stackrel{\d^\nds}{\longrightarrow} \ce_0(V)
$$ 
are complexes if and only if $F_{-}=0$.
In Riemannian signature each of these is an elliptic complex.
\end{proposition}

Evidently then we obtain detour complexes by composing the twisted de
Rham subcomplexes in the Proposition. For example if the connection
$\nds$ is anti-self-dual then there is a detour complex
\begin{equation}\label{asddet}
\ce^0(V)\stackrel{d^\nds}{\longrightarrow} \ce^1(V)\stackrel{2 \d^\nds_+d^\nds_+}{\longrightarrow} 
\ce_1(V)\stackrel{\d^\nds}{\longrightarrow}\ce_0(V).
\end{equation}
Similarly if $\nds$ is
instead self-dual then there is a detour complex
\begin{equation}\label{sddet}
\ce^0(V)\stackrel{d^\nds}{\longrightarrow} \ce^1(V)\stackrel{2 \d^\nds_-d^\nds_-}{\longrightarrow} 
\ce_1(V)\stackrel{\d^\nds}{\longrightarrow}\ce_0(V).
\end{equation}
The following result is a straightforward calculation.
\begin{proposition}\label{agree} 
The complexes \nn{asddet} and \nn{sddet} are special cases of the
twisted de Rham detour complex \nn{detseq} of Theorem \ref{twistthm}.
\end{proposition}

\section{Translation via the Yang-Mills detour complex}\label{detdiag}

 We may use the Theorem \ref{twistthm} to construct more exotic
differential complexes. The ideas here are partly inspired by
Eastwood's curved translation principle \cite{ER,Esrniconf} which
in turn is a geometric adaptation of the Jantzen-Zuckermann translation
functor from representation theory.

Consider the following general situation. Suppose that there are vector
bundles (or rather section spaces thereof) $\cB^0$, $\cB^1$, $\cB_1$ and
$\cB_0$ and differential operators $L_0$, $L_1$, $L^1$, $L^0$, $\cD$
and $\overline{\cD}$ which act as indicated in the following diagram:
 \newcommand{\fbbbb} {\mbox{$
\begin{picture}(40,2)(1.6,-2)
\put(1,1.0){\line(1,0){58}}
\end{picture}$}} 
$$
\begin{picture}(350,100)(12,-80)

\put(10,-25){[D]}

\put(70,0) {$ \ce^0(V)\hspace*{2pt}\stackrel{d^\nds}{\longrightarrow}
\hspace*{8pt} \ce^1(V) \stackrel{M^\nds}{\fbbbb\longrightarrow}
\hspace*{8pt} 
\ce_1(V)\hspace*{2pt} \stackrel{\d^\nds}{\longrightarrow}\hspace*{8pt} 
\ce_0(V)
$}                                                         

\put(77,-38){\vector(0,1){31}}                               
\put(66,-25){\scriptsize$L_0$}

\put(143,-38){\vector(0,1){31}}                               
\put(132,-25){\scriptsize$L_1$}

\put(249,-6){\vector(0,-1){31}}                              
\put(238,-25){\scriptsize$L^1$}

\put(315,-6){\vector(0,-1){31}}                              
\put(304,-25){\scriptsize$L^0$}

\put(70,-50)
{$
\hspace*{2pt} \cB^0{\phantom{(V)}} \stackrel{\cD}{\longrightarrow} 
\hspace*{8pt} \cB^1{\phantom{(V)}} \hspace*{-2pt}
\stackrel{M^{\cB}}{\fbbbb\longrightarrow} \hspace*{8pt}
\cB_1{\phantom{(V)}}
\stackrel{\overline{\cD}}{\longrightarrow}\hspace*{8pt}  \cB_0{\phantom{(V)}}
$}                                                         

\end{picture}
$$

\vspace{-8mm}

\noindent The top sequence is \nn{detseq} for a connection $D$ with curvature
$F$ and the operator $M^\cB:\cB^1\to \cB_1$ is defined to be the
composition $L^1M^{\nds}L_1$.  Suppose that the squares at each end
commute, in the sense that as operators $\cB^0\to \ce^1(V)$ we have
$d^{\nds}L_0=L_1\cD$ and as operators $\ce_1(V)\to \cB_0$ we have
$L^0\d^\nds=\overline{\cD}L^1$.  Then on $\cB^0$ we have
$$
M^{\cB} \cD = L^1 M^D L_1 \cD = L^1 M^D d^D L_0 = L^1\e(\d^D F) L_0 ,
$$
and similarly $\overline{\cD} M^{\cB}=- L^0 \i(\d^D F) L_0 $. Thus if
$D$ is Yang-Mills then the lower sequence, viz.
\begin{equation}\label{transcx}
\cB^0 \stackrel{\cD}{\longrightarrow}  \cB^1
\stackrel{M^{\cB}}{\longrightarrow} 
\cB_1
\stackrel{\overline{\cD}}{\longrightarrow}  \cB_0 ~,
\end{equation}
is a complex. 

{\bf Remarks:} Note that if the connection $D$ preserves a
Hermitian or metric structure on $V$ then we need only the single commuting 
square $d^{\nds}L_0=L_1\cD$ on $\cB^0$ to obtain such a complex; by taking formal adjoints we obtain a 
second commuting square $(L^0\d^\nds=\overline{\cD}L^1): \cB_1\to \cB_0$ where 
$\cB_0$ and $\cB_1$ are appropriate density twistings of the bundles dual to 
$\cB^0$ and $\cB^1$ respectively.

Obviously for \nn{transcx} to be a complex, it
is sufficient (and necessary) for $L^1M^{\nds}(d^\nds L_0-L_1 \cD)$ to
vanish on $\cB^0$ and for $(\overline{\cD} L^1-L^0\d^\nds) M^\nds L_1$
to vanish on $\cB^1$. 
\quad \endrk 

\subsection{The complex for (almost) Einstein scales} \label{einsect}

We work in the setting of conformal $n$-manifolds, $n\geq 3$. We will
construct here a diagram of the form [D] via the normal conformal
tractor connection. The standard tractor bundle is vector bundle with
a conformally invariant connection that we may view as arising as
an induced structure from the Cartan bundle and connection of
\cite{Cartan}.  In fact the Cartan connection is readily recovered
from the tractor connection, see \cite{CapGotrans} where such
connections and related calculus are described for the class of
parabolic geometries (which also includes, for example, CR geometry,
quaternionic structures and projective geometry).  For our current
construction it is not the normality of the tractor connection, in the
sense of \cite{CapGotrans,Cartan}, that is important. Rather the key
point is that it arises from a prolongation (as observed in
\cite{BEGo}) of a certain (finite type) partial differential operator
$P$ that we may take as the operator $\cD$ for the diagram [D]: In
terms of a metric $g$, this operator $P$ is given by
\begin{equation}\label{sc}
 P\sigma ={\rm TF}(\nd_{a} \nd_{b}\si+\Rho_{ab}\si ),
\end{equation}
where $\si\in \ce[1]$.  Modulo the trace part, this is the
differential operator which controls the conformal transformation of
the Schouten tensor.  In particular a metric $\si^{-2}\bg$ is Einstein
if and only if the scale $\si\in \ce[1]$ is non-vanishing and
satisfies $P \si =0$.  In order to be explicit we give a construction
of the tractor connection here, as it is the key to obtaining the
required commutative diagram. For further details see
\cite{CapGoluminy}. 

We write $J^k E[1]$ for the bundle of k-jets of germs of sections of
$E[1]$.  Considering, at each point of the manifold, sections which
vanish to first order at the given  point reveals a
canonical sequence,
$$
0\to S^2 T^*M\otimes E[1] \to J^2 E[1]\to J^1E[1]\to 0~.
$$
This is the jet exact sequence at 2-jets. Via the conformal metric
$\bg$, on a conformal manifold the bundle of symmetric covariant
2-tensors $S^2 T^*M$ decomposes directly into the trace-free part,
which we will denote $E^{1,1}$, and a pure trace part isomorphic to
$E[-2]$, that is $S^2 T^*M[1] = E^{1,1}[1]\oplus E[-1]$.  The {\em standard
  tractor bundle} $\bT$ may defined as the quotient of $J^2 E[1]$ by
the image of $E^{1,1}[1]$ in $J^2 E[1]$. By construction this is
invariant, it depends only on the conformal structure. Also by construction, 
it is an
extension of the 1-jet bundle
$$
0\to E[-1]\to \bT \to J^1 E[1]\to 0.
$$
 Note that there is a tautological operator $\bD:\ce[1]\to
\ce^0(\bT)$ which is simply the composition of the universal 2-jet
differential operator $j^2:\ce[1]\to \ce^0(J^2 E[1])$ followed by the canonical
projection $\ce^0(J^2 E[1])\to \ce^0(\bT)$. By construction this is
invariant. 

 Via a choice of metric $g$, and the Levi-Civita
connection it determines, we obtain a differential operator $\ce[1]\to
\ce[1]\oplus \ce^1[1]\oplus \ce[-1]$ by $\si\mapsto ( \si, \nd_a \si,
- \frac{1}{n}(\Delta + \J) \si )$ and this obviously determines an
isomorphism
\begin{equation}\label{split}
\ce^0(\bT) \stackrel{g}{\cong} \ce[1]\oplus \ce^1[1]\oplus
\ce[-1] ~.
\end{equation}
Changing to a conformally related metric $\widehat{g}=e^{2\om}g$
($\om$ a smooth function) gives a different isomorphism, which is
related to the previous by the transformation formula
\begin{equation}\label{transf}
\widehat{(\si,\mu_b,\tau)}=(\si,\mu_b-\si\Up_b,\tau+\bg^{bc}\Up_b\mu_c-
\tfrac{1}{2}\si\bg^{bc}\Up_b\Up_c),  
\end{equation}
where $\Upsilon:=d\om$. 
Now we define a connection on $\ce[1]\oplus \ce^1[1]\oplus
\ce[-1] $ by the formula
\begin{equation}\label{trconn}
\nd_a
\left( \begin{array}{c}
\si\\\mu_b\\ \rho
\end{array} \right) : =
\left( \begin{array}{c}
 \nabla_a \si-\mu_a \\
 \nabla_a \mu_b+ g_{ab} \rho +\Rho_{ab}\si \\
 \nabla_a \rho - \Rho_{ab}\mu^b  \end{array} \right) 
\end{equation}
where, on the right-hand-side $\nd$ is the Levi-Civita connection for
$g$.  Obviously this determines a connection on $\bT$ via the
isomorphism \nn{split}. What is more surprising is that if we repeat
this using the metric $\widehat{g}$, conformally related to $g$, in
\nn{split} and \nn{trconn} we obtain the {\em same} connection on
$\bT$. This may easily be verified by, for example, directly
calculating that under such a conformal change the right-hand side of
\nn{trconn} transforms in exactly the same way as a (1-form valued)
invariant section of $\bT$. That is it transforms according to \nn{transf}.
The canonical connection on $\bT$, so constructed, depends only on the
conformal structure and is known as the {\em (normal standard) tractor
connection}. 
In what
follows we will use \nn{split} without further explicit comment.
There is also a conformally invariant {\em tractor metric} $\h$ on
$\bT$ given (as a quadratic form) by $(\si,~\mu,~\rho)\mapsto
\bg^{-1}(\mu,\mu)+2\si \rho$. This is preserved by the connection and
has signature $(p+1,q+1)$ (corresponding to $\bg$ of signature
$(p,q)$).  

Note that, given a metric $g$, through \nn{split} the tautological
invariant operator $\bD$ from above is given by the explicit formula
$$ \ce^0[1]\to \ce^0(\bT) \quad \quad \quad \si\mapsto (\si,\nd_a \si
, -\frac{1}{n}(\Delta \si +\J \si)) .
$$ This is a called a differential splitting operator since through
the jet projections there is conformally invariant surjection $X:
\ce(\bT)\to \ce[1]$ which inverts $\bD$. There is also a
differential splitting operator 
$$
 E:\ce^{1,1}[1] \to \ce^1(\bT) \quad \quad \quad \psi_{ab} \mapsto
(0,\psi_{ab},-(n-1)^{-1}\nd^b\psi_{ab})
$$ (cf.\ \cite{Esrniconf}). An easy calculation verifies that this 
also is conformally invariant. 
We have the following. 
\begin{proposition}\label{eincomm}
With $\nd$ denoting the tractor connection on $\ce^0(\bT)$ we have 
$$
\nd \bD = E P   ~.
$$ as differential operators on $\ce[1]$. For $\si\in \ce[1]$, $\bD\si$
is parallel if and only if $P\si=0$. 
\end{proposition}
\noindent{\bf Proof:} The second statement is immediate from the
first.  A straightforward calculation verifies that either composition
applied to $\si\in \ce[1]$ yields
$$
\left( \begin{array}{c}
0\\ {\rm TF}(\nd_{a} \nd_{b}\si+\Rho_{ab}\si ) \\ 
-\frac{1}{n}\nd_a(\Delta \si +J \si)-P_a{}^c\nd_c\si
\end{array} \right) 
$$ \quad $\blacksquare$ \\ 
In fact if a section $I\in \ce^0(\bT)$ is parallel then $I=\bD\si$ for
some $\si\in \ce[1]$ so a conformal manifold with a parallel tractor
is {\em almost Einstein} in the sense that it has a section of
$\ce[1]$ that gives an  Einstein scale on an open dense subset (see
\cite{goalmost} for further details).

Since the tractor connection is orthogonal (for the conformally
invariant tractor metric $h$ given above) the formal adjoints of the 
operators above give another commutative square of operators. That is with 
$$\begin{array}{ccll}
 && P^*: \ce_{1,1}[-1]\to \ce_0[-1] & 
\phi_{ab}\mapsto \nd^{a}\nd^b \phi_{ab}+ \Rho^{ab}\phi_{ab}  \\ 
&& E^*: \ce_1(\bT) \to \ce_{1,1}[-1] & (\alpha_a,~\nu_{ab},~\tau_a) \mapsto \nu_{(ab)_0}
+\frac{1}{n-1}\nd_{(a}\alpha_{b)_0} \\
 && \bD^*: \ce_0(\bT) \to \ce_0[-1]& (\si,~\mu_b,~\rho) \mapsto \rho -\nd^a\mu_a-\frac{1}{n}(\Delta \si+ J\si)\\
&&  \d^\nd: \ce_{1}(\bT)\to \ce_0(\bT) & \Phi_{aB}\mapsto -\nd^a \Phi_{aB}~,
\end{array}
$$ 
\noindent where $\ce_{1,1}$ denotes the space of sections of $E^{1,1}\otimes E[4-n]$,   
we have $\bD^*\delta^\nd= P^* E^*$ on $\ce_{1}(\bT)$.

Finally observe that 
the curvature of the tractor connection, as calculated directly 
from \nn{trconn}, is
$$
\Om_{ab}{}^C{}_D=
\left(\begin{array}{lll}
0&0&0\\
A^{c}{}_{ab} & C_{ab}{}^c{}_d & 0\\
0 & -A_{dab} & 0
\end{array}\right)
$$
and hence (see e.g.\ \cite{GoNur} for further details),  
\begin{equation}\label{divtr}
\nd^a \Om_{ab}{}^C{}_D= \left(\begin{array}{lll}
0&0&0\\
B^{c}{}_{b} & (n-4)A_{b}{}^c{}_d & 0\\
0 & -B_{db} & 0
\end{array}\right)
\end{equation}
 where, on the left-hand side, $\nd$ is the Levi-Civita connection
coupled with the tractor connection on $\End(\bT)$ induced from
\nn{trconn}.  Let us say that a {\em pseudo-Riemannian} manifold is
{\em semi-harmonic} 
 if its tractor curvature is Yang-Mills, that
is $ \nd^a \Om_{ab}{}^C{}_D=0$. Note that in dimensions $n\neq 4$ this
is not a conformally invariant condition and a semi-harmonic
space is a Cotton space that is also Bach-flat.   From our 
observations above, the
semi-harmonic condition is conformally invariant in dimension 4
and according to the last display we have the following
result. 
\begin{lemma}\label{trYM}
 In dimension 4 the tractor connection \nn{trconn} is a Yang-Mills
connection if and only if the structure is Bach-flat.
\end{lemma}
\noindent This result is not new and equivalent observations have been
known in the literature for some time \cite{Merk,MB,KL}. 
It brings us to 
 the following.  Let us write $M^{\bT}$ for the composition $E^* M^\nd E $.
On $h\in \ce^{1,1}[1]$ we have
$$
(M^{\bT} h)_{ab}
= -TFS\big(\nd^c(\nd_c h_{ab}-\nd_a h_{cb}) 
-\frac{1}{n-1}\nd_a \nd^c h_{b c} +
 C_a{}^c{}_b{}^d h_{cd} 
\big ),
$$
and the following results. 
\begin{theorem}\label{paracase}
The sequence
\begin{equation}\label{Einseq}
 \ce^0[1]\stackrel{P}{\to} \ce^{1,1}[1] 
\stackrel{M^{\bT}}{\longrightarrow} \ce_{1,1}[-1]
\stackrel{P^*}{\to} \ce_0[-1]
\end{equation}
has the following properties. \\
\IT{i} It is a formally self-adjoint sequence of
differential operators and, 
for $\si\in \ce^0[1]$ 
\begin{equation}\label{MP}
(M^{\bT} P \si)_{ab} =
-TFS\big( B_{ab}\si - (n-4) A_{abc}\nd^c\si \big) ,
\end{equation}
where $TFS(\cdots)$ indicates the trace-free symmetric part of the
tensor concerned.  In particular it is a complex on
semi-harmonic manifolds. \\ 
\IT{ii} In the case of Riemannian
signature the complex is elliptic. \\ 
\IT{iii} In dimension 4,
\nn{Einseq} is a sequence of conformally invariant operators and it is a
complex if and only if the conformal structure is Bach-flat.
\end{theorem}
\noindent{\bf Proof:}
Setting 
$$
\cD=P, \quad L_0=\bD, \quad L_1=E 
$$ we have the situation of the translation diagram [D] above, with 
the right square given by formal
adjoints of these operators, and the tractor bundle connection pair
$(\bT,\nd)$ used for $(V,D)$ in the top row. That is:

{\begin{picture}(350,100)(12,-80)

\put(70,0) {$ \ce^0(\bT)\hspace*{2pt}\stackrel{d^\nd}{\longrightarrow}
\hspace*{8pt} \ce^1(\bT) \stackrel{M^\nd}{\fbbbb\longrightarrow}
\hspace*{8pt} 
\ce_1(\bT)\hspace*{9pt} \stackrel{\d^\nd}{\longrightarrow}\hspace*{11pt} 
\ce_0(\bT)
$}                                                         

\put(77,-38){\vector(0,1){31}}                               
\put(66,-25){\scriptsize$\bD$}

\put(143,-38){\vector(0,1){31}}                               
\put(132,-25){\scriptsize$E$}

\put(260,-6){\vector(0,-1){31}}                              

\put(248,-25){\scriptsize$E^*$}
\put(335,-6){\vector(0,-1){31}}                              
\put(324,-25){\scriptsize$\bD^*$}

\put(70,-50)
{$
\hspace*{2pt} \ce^0[1]  \stackrel{P}{\longrightarrow} 
\hspace*{8pt} \ce^{1,1}[1]  \hspace*{2pt}
\stackrel{M^{\bT} }{\fbbbb\longrightarrow} \hspace*{8pt}
 \ce_{1,1}[-1] 
\stackrel{P^*}{\longrightarrow}\hspace*{8pt} \ce_0[-1] 
$}                                                         

\end{picture}}

\vspace{-8mm}

\noindent 
By construction the lower sequence \nn{Einseq} is formally
self-adjoint, and in dimension 4 conformally invariant.  If the
structure is semi-harmonic then the upper sequence is a complex and
hence, from the commutativity of the diagram, \nn{Einseq} is a
complex. In particular on Bach-flat 4-manifolds we obtain a complex.
On the other hand, from \nn{MP} it follows that in dimension 4 we
obtain a complex only if the structure is Bach-flat.

From \nn{trconn} we calculate 
$d^\nabla$ on the range of $E$ to obtain  
\begin{equation}\label{twD}
d^\nabla \left(\begin{array}{c} 0\\
\nu\\
\tau\end{array}\right) = \left(\begin{array}{c} 0\\
Q \nu \\
 \ast \end{array}\right) 
\end{equation}
 where, for $\nu\in\ce^{1,1}[1]$, we have
$\tau=-\frac{1}{(n-1)}\nd^b\nu_{ab}$,
$Q$ is given by
$$
(Q\nu)_{abc}=2\nd_{[a}\nu_{b]c}+2\bg_{c[a}\tau^{\phantom{a}}_{b]} ,
$$ and we do not need the details of the term indicated by $\ast$.  
It follows immediately from \nn{twD}, and the formulae for the tractor
metric, that we have $M^{\bT}:= E^* M^\nabla E= Q^* Q+\LOT$ where $Q^*$
denotes the formal adjoint of the operator $Q$. In Riemannian signature 
the leading symbol of
$Q^* Q$ has the same kernel as the leading symbol of $Q$, and it
follows easily that the complex is elliptic. The ``lower order terms''
(indicated by $\LOT$) in $M^{\bT}$ arise simply from the tractor
curvature in the formula for $M^D$ and amount to an action by the Weyl
curvature. Including this yields the explicit formula for $M^{\bT}$
given above the Theorem. 
The expression
\nn{MP} for the composition $M^{\mathbb{T}} P$ follows from this by a short
direct calculation. (The calculation is even
simpler if the result of Lemma \ref{algact} is imported).
\quad $\blacksquare$

\vspace{2mm}

\noindent{\bf Remark:} 
Note that the formula for $M^{\bT}$ is
 closely related to, but not the same as, the operator which
arises from deformations of Einstein structures (for the latter see
e.g.\ \cite{Besse} and references therein). It should be valuable to
expose the geometric meaning of the first cohomology of the
sequence \nn{Einseq}.
\quad \endrk

\begin{corollary}\label{EimpliesB}
  Einstein 4-manifolds are Bach-flat. 
 \end{corollary}
\noindent{\bf Proof:} If a non-vanishing density $\si$ is an
Einstein scale then, calculating in that scale, we have $M^{\bT} P\si
= -B \si$, where $B$ is the Bach tensor. On the other hand if $\si$ is
an Einstein scale then $P \si =0$ (see \nn{sc}).  \quad $\blacksquare$

\noindent{\bf Remarks:} In fact, more generally, almost Einstein
manifolds are also necessarily Bach flat. Since an almost Einstein
manifold has an Einstein scale on an open dense subspace, this follows
by continuity of the Bach tensor. (The higher dimensional extension
of this result is that even dimensional almost Einstein manifolds have
vanishing Fefferman-Graham obstruction tensor, see e.g.\
\cite{goalmost} and references therein.)  The result that Einstein
metrics are Bach-flat is well-known by other means (see e.g.\
\cite{KNT,GoNur}). Nevertheless we feel the detour complex gives an
interesting route to this.  In any dimension Einstein metrics satisfy
$P_{ab}=\frac{1}{n}J g_{ab}$ with $J$ constant, so it follows from the
definitions of the Cotton tensor and the Bach tensor \nn{bach} that
Einstein metrics are semi-harmonic. Thus there are many examples of
semi-harmonic manifolds.

\subsection{The twistor spinor complex}\label{tws}

We assume here that we have a conformal spin structure. This is no
 restriction locally.  For the purpose of being self-contained and
 having the results in a uniform notation we derive the
 basic spinor identities we require. An alternative treatment
 for many of these may be found in, e.g.\
 \cite{Fried-book}. We will use the spin-tractor connection
 below. This is often termed the local twistor connection
 \cite{otnt,BFGK}. The notation we use (and the basic tractor tools)
 follows \cite{TomLeb} which presents a spin-tractor calculus
 developed by the first author and Branson. Following that source we
 write $\S$ for the basic spinor bundle and $\overline{\S}={\S}[-n]$
 (i.e.\ the bundle that pairs globally in an invariant way with $\S$
 on conformal $n$-manifolds). Evidently the weight conventions here
 give $\S$ a ``neutral weight''. In terms of, for example, the Penrose
 weight conventions
 $\S=E^\lambda[-\frac{1}{2}]=E_{\lambda}[\frac{1}{2}]$, where
 $E^\lambda$ denotes the basic contravariant spinor bundle in
 \cite{otnt}.

We write ${\rm Tw}$ for the so-called twistor bundle, that is the
 subbundle of $T^*M\otimes \S[1/2]$ consisting of form spinors $u_a$
 such that $\gamma^a u_a=0$, where $\gamma_a$ is the usual Clifford
 symbol.  We use $\S$ and ${\rm Tw}$ also for the section spaces of
 these bundles.
The {\em twistor operator} is the conformally invariant
 Stein-Weiss gradient
$$
{\bf T}: \S[1/2]\to {\rm Tw}
$$
given explicitly by
$$
\psi\mapsto \nd_a \psi +\frac{1}{n}\gamma_a\gamma^b \nd_b \psi ~.
$$ 
The main result of this section is that this completes to a
differential complex as follows.
\begin{theorem}\label{twcase}
On semi-harmonic pseudo-Riemannian $n$-manifolds $n\geq 4$ we have a
differential complex
\begin{equation}\label{spcx}
 \S[1/2]  \stackrel{\bf T}{\to} {\rm Tw} 
\stackrel{M^\Sigma}{\longrightarrow} \overline{\rm Tw}
\stackrel{{\bf T}^*}{\to} \overline{\S}[-1/2] ,
\end{equation} 
 where ${\bf T}$ is the usual twistor operator, ${\bf T}^*$ its
formal adjoint, and ${M^\Sigma}$ is a third order operator and given by the formula
\nn{Neq} below. The sequence is formally self-adjoint and in the case
of Riemannian signature the complex is elliptic.

In dimension 4 the sequence \nn{spcx} is conformally invariant and it
is a complex if and only if the conformal structure is Bach-flat.
\end{theorem}
\noindent{\bf Remarks:} Of course on a fixed pseudo-Riemannian
manifold we may ignore the conformal weights. 

Note also that, for
example in dimension 4, under the chirality decomposition of this
sequence, we get the two complexes
$$
\S_\pm[1/2]  \stackrel{\bf T}{\to} {\rm Tw}_\pm 
\stackrel{M^\Sigma}{\longrightarrow} \overline{\rm Tw}_\mp
\stackrel{{\bf T}^*}{\to} \overline{\S}_\mp[-1/2] ,
$$ by the restriction of the operators ${\bf T}$, ${\bf T}^*$,
and ${M^\Sigma}$.

If we were to apply the construction below in dimension 3 then we
would obtain a trivial operator ${M^\Sigma}$. In this dimension  the
BGG sequence (see section \ref{back}) takes the form
\nn{spcx}, where the middle operator is of second order. 
  \quad \endrk

In the calculations which follow it will often be convenient to use
abstract indices for the form bundles while at the same time not using
any indices for the spinor bundles.  We have already done this
implicitly above, for example in the  formula for the twistor operator
which, in this notation, acts ${\bf T}_a: \S[1/2]\to {\rm Tw}_a$.  From
the usual gamma matrices $\gamma^a$ satisfying
$$
\gamma^a \gamma^b+\gamma^b \gamma^a =-2 \bg^{ab} {\rm Id} 
$$
we switch to the  symbols $\beta:=\gamma/\sqrt{2}$, so that 
\begin{equation}\label{crel}
\beta^a \beta^b+\beta^b \beta^a =- \bg^{ab} {\rm Id} ,
\end{equation}
this simplifies certain formulae in the following discussion. 
 We denote the corresponding  Dirac operator by $\D:=\beta^a\nd_a$.

Given a metric $g$ from the conformal class the spin-tractor bundle
$\Sigma$  is given by
$$
\Sigma \stackrel{g}{\cong} \S [1/2] \oplus \S [-1/2]
$$ 
where $\S$.  In the conformally
 related metric $\widehat{g}=e^{2\om}g$ we have a similar isomorphism
 and
\begin{equation}\label{spintrtr}
\left(\begin{array}{c}\widehat{\psi} \\ \widehat{\phi} \end{array} \right) = \left(\begin{array}{c}\psi\\\phi + \Upsilon_c\beta^c \psi \end{array} \right)
\end{equation}
where $\Upsilon=d\om$. In terms of the $g$-splitting  
the normal conformal spin-tractor connection is given by 
\begin{equation}\label{spinconn}
 \nd_a \left(\begin{array}{c}\psi\\\phi \end{array} \right)
=\left(\begin{array}{c}\nd_a\psi+\beta_a \phi\\
\nd_a\phi+P_{ab}\beta^b\psi \end{array} \right)~.
\end{equation}
On the right side $\nd$
means the usual Levi-Civita (spin) connection, while on the left the
same notation is used for the spin-tractor connection.
It is an easy exercise to verify directly that this is a conformally
invariant connection. The normality follows from the characterisation
of normal tractor connections (for irreducible parabolic geometries)
given in Theorem 1.3 of \cite{CapGoluminy}.  

The invariant pseudo-Hermitian form on spin-tractors is given by
$$
\langle \phi, \underline{\psi} \rangle+\langle \psi,\underline{\phi}\rangle
$$ for a pair spin-tractors $(\phi,\psi)$,
$(\underline{\phi},\underline{\psi})$, and where $\langle\cdot,\cdot
\rangle$ is the usual Hermitian form on spinors (which is compatible
with Clifford multiplication and is preserved by the Levi-Civita spin
connection). It is readily verified that this is invariant under the
transformations \nn{spintrtr} and that it is preserved by the
spin-tractor connection \nn{spinconn}.  We subsequently calculate in a
metric scale $g$ without further comment.

We construct two differential splitting operators: 
$$
L_0: \S[1/2]\to \ce^0(\Sigma)
$$
is given by
\begin{equation}\label{L0}
\psi\mapsto \left(\begin{array}{c}\psi\\ \frac{2}{n}\D\psi \end{array} \right)~;
\end{equation}
$$
L_1: {\rm Tw}\to \ce^1(\Sigma)
$$
is given by
\begin{equation}\label{L1}
\psi_a\mapsto \left(\begin{array}{c}\psi_a\\ \frac{2}{n-2}(\D\psi_a-\frac{1}{n-1}\beta_a\nd^b\psi_b) \end{array} \right)~.
\end{equation}
It is a straightforward exercise to verify that these transform
according to \nn{spintrtr} and so are expressions, in the metric scale
$g$, for conformally invariant operators.
An essential feature of these operators is the following commutativity result. 
\begin{proposition}\label{spincomm}
With $\nd$ denoting the spin-tractor connection on $\ce^0(\Sigma)$ we have 
$$
d^\nd L_0 = L_1 {\bf T}   ~,
$$ as differential operators on $\S [1/2]$. For $\psi\in \S [1/2]$, $L_0 \psi$
is parallel if and only if $ {\bf T} \psi =0$.
\end{proposition}
\noindent The last comment follows immediately from the commutativity
of the square given that $L_1$ is a differential splitting operator
(and so, in particular, $L_1 \psi_a=0 \Rightarrow \psi_a =0$). A
correspondence between parallel spin-tractors and twistor spinors
dates back to \cite{Fried89}. The extra information in the proposition
is that the operator $L_1$ is a conformally invariant tractor
splitting operator.  We shall postpone the proof of Proposition
\ref{spincomm}, as we prefer to first complete the proof of Theorem
\ref{twcase}.

Suppose that for a metric $g$, from the
 conformal class, the tractor connection is semi-harmonic. Recall
 this is exactly the condition that the normal tractor connection is
 Yang-Mills.  It follows immediately that the spin-tractor connection
 is also Yang-Mills, since this is induced from the same principal
 connection simply pulled back to the 2-1 covering Spin$(p+1,q+1)$-principal
 bundle.  (Equivalently they arise from the same Cartan connection,
 this is the usual picture \cite{CapGotrans} and sufficient to see
 this result. From the Cartan picture one may easily extend to a
 principal bundle and connection from which the tractors are induced,
 and this is simply an alternative framework.)  As observed above, in
 dimension 4 the Yang-Mills condition is exactly the condition that
 the metric (or conformal structure) is Bach-flat. Thus, from the
 Proposition, the first part of Theorem \ref{twcase} follows immediately from
 the commutative diagram below where
\begin{equation}\label{Neq}
{M^\Sigma}= L^1 M^\nd L_1 .
\end{equation}
{\begin{picture}(350,100)(12,-80)

\put(70,0) {$ \ce^0(\Sigma)\hspace*{2pt}\stackrel{d^\nd}{\longrightarrow}
\hspace*{8pt} \ce^1(\Sigma) \stackrel{M^\nd}{\fbbbb\longrightarrow}
\hspace*{8pt} 
\ce_1(\Sigma)\hspace*{2pt} \stackrel{\d^\nd}{\longrightarrow}\hspace*{8pt} 
\ce_0(\Sigma)
$}                                                         

\put(77,-38){\vector(0,1){31}}                               
\put(66,-25){\scriptsize$L_0$}

\put(143,-38){\vector(0,1){31}}                               
\put(132,-25){\scriptsize$L_1$}

\put(264,-6){\vector(0,-1){31}}                              
\put(253,-25){\scriptsize$L^1$}

\put(335,-6){\vector(0,-1){31}}                              
\put(324,-25){\scriptsize$L^0$}

\put(70,-50)
{$
\hspace*{2pt} \S[1/2]  \stackrel{{\bf T}}{\longrightarrow} 
\hspace*{8pt} {\rm Tw} {\phantom{(\Sigma)}} \hspace*{-2pt}
\stackrel{{M^\Sigma} }{\fbbbb\longrightarrow} \hspace*{8pt}
 \overline{\rm Tw}  {\phantom{(\Sigma)}}
\stackrel{{\bf T}^*}{\longrightarrow}\hspace*{8pt} \overline{\S}[-1/2] {\phantom{(\Sigma)}}
$}                                                         

\end{picture}}

\vspace{-8mm}

\noindent 
Using that $\Sigma$ is a self-dual bundle, the operators in the square at
the right end of the diagram are defined as the formal adjoints of the
operators in the first square. So all squares commute and both
horizontal sequences are formally self-adjoint.

\newcommand{\K}{\mbox{\bf K}} 

To establish ellipticity we require the leading term of the operator
${M^\Sigma}$.  Applying the spin-tractor twisted exterior derivative to
$L_1 t$, for $t\in {\rm Tw}$, we obtain a result of the form
\begin{equation}\label{4ell}
\left(\begin{array}{c}\K t\\
\ell \D \K t + m \beta \d^\nabla \K t + curvature \end{array}\right)
\end{equation}
with $\K t$ given by 
$$
2\nabla_{a_1} t_{a_2}+ \frac{4}{n-2} \beta_{a_1}(\D t_{a_2}-\frac{1}{n-1}\beta_{a_2}\nabla^ct_c),
$$ where the indices $a_1a_2$ are implicitly skewed over, $\ell$ and
$m$ are constants, and $\d^\nabla$ is the spin-Levi-Civita connection
twisted interior derivative. By construction $\K$ is an invariant
operator $\K:{\rm Tw}\to {\rm Tw}^2$, where ${\rm Tw}^2$ is the
subbundle of $E^2\otimes \S[1/2]$ consisting of spin-forms annihilated
by interior multiplication by $\beta$.  It is a straightforward
exercise (or one may use the BGG machinery of \cite{CSSannals}) to construct a differential splitting operator $L_2:{\rm
Tw}^2\to E^2\otimes \Sigma $. This has the form $s\mapsto (s~,~\ell \D
s + m \beta\d^\nabla s)$ (cf. $L_1$) where $\ell\neq 0$. On the other
hand in the conformally flat case it follows easily, from the
uniqueness of conformal differential operators, that $L_2 \K= d^\nabla
L_1$. Thus we obtain the form of the bottom slot of \nn{4ell}.

Since the leading term of
${M^\Sigma}$ is obtained by composing $d^\nabla L_1$ with its formal
adjoint, and $\D$ is formally self-adjoint, it follows that the
leading term of the operator ${M^\Sigma}$ is of the form $\K^* \pi_2\D
\K$. Here $\K^*$ denotes the formal adjoint of $\K$ and the projection
$\pi_2$ is the projection of spinor-valued two-forms to the space $Tw^2.$

The operator $R_2=\pi_2 D$ on $Tw^2$ is an elliptic and self adjoint
operator,
it is a higher spin analogue of the Dirac operator. Similarly,
if $\pi_1$ denotes the projection of spinor-valued one-forms to $Tw,$
the operator $\R_1=\pi_1 D$ is an elliptic self adjoint operator on $Tw$
(usually called the Rarita-Schwinger operator).
Moreover, $R_1 \K^*$ is a multiple of $\K^*R_2.$
Hence the leading term of the operator ${M^\Sigma}$ is a multiple
of $R_1\K^* \K.$
                                                                                
Since elements of ${\rm Tw}^2$ are annihilated by interior Clifford
 multiplication, it follows from the formula for $\K$ that  the symbol
$\si_\xi(\K^*)$ is simply interior multiplication by $\xi$, $\i(\xi)$.

Without loss of generality we may suppose that $\xi$ is a unit
vector.
It is well known that the Rarita-Schwinger operator (on the flat space) can
be composed with a third order constant coefficient operator
to give the square of the Laplace operator. Hence the symbol
of $\R_1$ can be multiplied on the left to get $|\xi|^4=1.$
By this left multiplication,
$\si_\xi (\K^* R_2 \K)t=\si_\xi(R_1 \K^* \K)t=0 $ implies  $\i(\xi)\sigma_\xi
(\K)t
 =0$. Now the explicit formula for $ \sigma_\xi (\K)t $ is
$$
\e(\xi) t+\frac{2}{n-2}\e(\beta)(\beta_\xi
t-\frac{1}{n-1}\e(\beta)\i(\xi)t),
$$
where $\e(\cdot)$ indicates exterior multiplication.

Contracting with $\xi$ and setting to zero we obtain ($(n-3)$ times)
$$
(n-1)t=n\big(\e(\xi)\i(\xi)t+\frac{2}{n}\e(\beta)\beta_\xi
\i(\xi)t\big) ~.
$$
The right hand side here
is a multiple of $ \sigma_\xi ({\bf T}) \i(\xi) t$. Thus
 $t$ is in the range of $\si_\xi ({\bf T})$, as required.

\vspace{2mm}

\noindent{\bf Completing the proof and remarks:} From \nn{divtr},
\nn{L0}, and Lemma \nn{algact} it follows easily that in dimension 4
the composition $ {M^\Sigma}{\bf T}$ on $\phi$, a section of $\S[1/2]$,
is, up to a non-zero multiple, a Clifford multiplication of the Bach
tensor $B_{ab}\beta^b\phi$. Thus the formally self-adjoint sequence
\nn{twcase} is a complex if and only if the structure is Bach-flat, as
claimed in the theorem.  If $\phi$ is a twistor spinor (i.e.\ ${\bf
T}\phi=0 $) then this Clifford action of the Bach tensor on $\phi$
obviously vanishes. In Riemannian signature it is in fact
straightforward to recover a parallel standard tractor from the
parallel spin-tractor corresponding to a twistor spinor.  Thus
Riemannian manifolds admitting a twistor spinor are almost Einstein
and so Bach flat. In fact the last conclusion here is well-known \cite{BFGK}.
\quad \endrk

\vspace{2mm}

\noindent{\bf Proof of Proposition \ref{spincomm}:} Let $\psi$ be a
section of $\S[1/2]$.  Using the formula \nn{L0} for the splitting
operator and the expression \nn{spinconn} for the spin-tractor
connection we have
$$
\nd_a L_0 \psi= \nd_a \left(\begin{array}{c}\psi\\ \frac{2}{n}\D\psi \end{array} \right)
=\left(\begin{array}{c}\nd_a\psi+ \frac{2}{n}\beta_a \D\psi \\
\frac{2}{n} \nd_a\D \psi+P_{ab}\beta^b\psi \end{array} \right) 
$$ Recalling that $\gamma =\sqrt{2} \beta $ and that $\D=\beta^a\nd_a$,
we thus have
$$
\nd_a L_0 \psi
 =
 \left(\begin{array}{c} {\bf T}_a\psi \\
\frac{2}{n} \nd_a\D \psi+P_{ab}\beta^b\psi \end{array} \right) ~.
$$
From \nn{L1} it is clear that it remains to show that 
\begin{equation}\label{pl}
\frac{2}{n} \nd_a\D \psi+P_{ab}\beta^b\psi = \frac{2}{n-2}(\D {\bf T}_a\psi
-\frac{1}{n-1}\beta_a\nd^b{\bf T}_b\psi) .
\end{equation}

Let us note some simpler identities first. First for the Levi-Civita 
spin-connection the curvature on a spinor $\phi$ is given by 
$$
[\nd_a,\nd_b]\phi = -\frac{1}{2}R_{abcd}\beta^c\beta^d\phi ~,
$$
where $R$ is the usual Riemannian curvature. 
Then from the Bianchi identities and the Clifford relation \nn{crel} we get 
\newcommand{\Sc}{\rm Sc}
\begin{equation}\label{bcurv}
R_{abcd}\beta^b\beta^c\beta^d =\Ric_{ab}\beta^b \quad
R_{abcd}\beta^a\beta^b\beta^c\beta^d = -\Sc/2 
\end{equation}
Next an elementary calculation shows that ${\bf T}^*$, the 
formal adjoint of ${\bf T}$, is given on $Tw$ by
$
\phi_a\mapsto -\nd^a\phi_a.
$
Thus, for $\psi$ in $\S [1/2]$,
$$
- {\bf T}^* {\bf T} \psi =\Delta \psi +\frac{2}{n}\D^2\psi, 
$$ (where $\Delta:=\nd^a\nd_a$)  since of the spin-connection preserves the Clifford
symbols.  On the other hand since $\D^2=\beta^a\nd_a\beta^b\nd_b =
\beta^a\beta^b\nd_a\nd_b$ and the $\beta$s anti-commute up to a trace,
as in \nn{crel}, while the $\nd$s commute up to curvature we obtain
$$
\D^2\psi =-\frac{1}{2}\Delta \psi -\frac{1}{4} R_{abcd}\beta^a\beta^b\beta^c\beta^d  \psi
$$
and so using \nn{bcurv} we come to 
$$
\Delta \psi = - 2 \D^2\psi + \frac{1}{4}\Sc \cdot \psi .
$$
This with the expression above for ${\bf T}^* {\bf T} $ gives
\begin{equation} \label{keyid}
-{\bf T}^* {\bf T} \psi = 2\big(\frac{1-n}{n}\big) \D^2 \psi +\frac{1}{4}\Sc\cdot \psi ~.
\end{equation}

We are now ready to calculate the left-hand side of \nn{pl}. 
Applying $\beta^b\nd_b$ to the defining identity 
$\nd_a\psi =-\frac{2}{n}\beta_a\D \psi+{\bf T}_a\psi$ we get 
$$ \beta^b\nd_b \nd_a \psi =-\frac{2}{n} \beta^b\beta_a\nd_b \D \psi
+\beta^b\nd_b{\bf T}_a \psi .
$$
Commuting the derivatives on the left and writing $\D$ as a 
shorthand for $\beta^b \nd_b$ (applied to e.g.\ ${\bf T}\psi$), we obtain 
$$
\nd_a\D \psi +\frac{1}{2} R_{abcd}\beta^b\beta^c\beta^d \psi
= \frac{2}{n}\nd_a \D\psi +\frac{2}{n}\beta_a \D^2\psi +\D {\bf T}_a\psi ~.
$$
Next rearranging the  terms and using \nn{bcurv} gives
$$ \big(\frac{n-2}{n}\big)\nd_a\D \psi= \frac{2}{n} \beta_a \D^2\psi
-\frac{1}{2}\Ric_{ab}\beta^b\psi +\D {\bf T}_a\psi.
$$
Using now the identity \nn{keyid} from above to substitute for $\D^2 \psi$
yields 
$$ \big(\frac{n-2}{n}\big)\nd_a\D \psi=
-\frac{1}{2}\big(\Ric_{ab}-\frac{1}{2(n-1)}\Sc \bg_{ab} \big)\beta^b \psi
+\D{\bf T}_a\psi +\frac{1}{n-1}\beta_a{\bf T}^* {\bf T} \psi ~.
$$ But multiplying this through with $2/(n-2)$ and using that the
Schouten tensor $P_{ab}=
\frac{1}{n-2}(\Ric_{ab}-\frac{1}{2(n-1)}\Sc \bg_{ab}  )$, and once
again that ${\bf T}^*\phi=-\nd^a\phi_a$, this gives exactly the
expression \nn{pl}, which is thus seen to be an identity.  \quad
$\blacksquare$

\end{document}